# HIGH MOMENT PARTIAL SUM PROCESSES OF RESIDUALS IN GARCH MODELS AND THEIR APPLICATIONS[1]

By Reg Kulperger and Hao Yu

*University of Western Ontario*

In this paper we construct high moment partial sum processes based on residuals of a GARCH model when the mean is known to be 0. We consider partial sums of $k$th powers of residuals, CUSUM processes and self-normalized partial sum processes. The $k$th power partial sum process converges to a Brownian process plus a correction term, where the correction term depends on the $k$th moment $\mu_k$ of the innovation sequence. If $\mu_k = 0$, then the correction term is 0 and, thus, the $k$th power partial sum process converges weakly to the same Gaussian process as does the $k$th power partial sum of the i.i.d. innovations sequence. In particular, since $\mu_1 = 0$, this holds for the first moment partial sum process, but fails for the second moment partial sum process. We also consider the CUSUM and the self-normalized processes, that is, standardized by the residual sample variance. These behave as if the residuals were asymptotically i.i.d. We also study the joint distribution of the $k$th and $(k+1)$st self-normalized partial sum processes. Applications to change-point problems and goodness-of-fit are considered, in particular, CUSUM statistics for testing GARCH model structure change and the Jarque–Bera omnibus statistic for testing normality of the unobservable innovation distribution of a GARCH model. The use of residuals for constructing a kernel density function estimation of the innovation distribution is discussed.

**1. Introduction and results.** In nonlinear time series and in particular econometric and discrete time financial modeling, Engle's [12] ARCH model plays a fundamental role; see [10], or the volume edited by Rossi [20]. The ARCH model has been generalized to GARCH by Bollerslev [5].

Received April 2003; revised November 2004.

[1]Supported by grants from the Natural Sciences and Engineering Research Council of Canada.

A GARCH$(p, q)$ sequence $\{X_t, -\infty < t < \infty\}$ is of the form

$$(1.1) \qquad\qquad X_t = \sigma_t \varepsilon_t$$

and

$$(1.2) \qquad\qquad \sigma_t^2 = \alpha_0 + \sum_{i=1}^{p} \alpha_i X_{t-i}^2 + \sum_{j=1}^{q} \beta_j \sigma_{t-j}^2,$$

where

$$(1.3) \qquad \alpha_0 > 0, \qquad \alpha_i \geq 0, \qquad 1 \leq i \leq p, \qquad \beta_j \geq 0, \qquad 1 \leq j \leq q$$

are constants, and the innovations process $\{\varepsilon_t, -\infty < t < \infty\}$ is a sequence of i.i.d. random variables (r.v.'s). When $\varepsilon_0$ has a finite $k$th moment, denote $\mu_k = E(\varepsilon_0^k)$. A usual GARCH model assumption also is:

$$(1.4) \qquad \begin{array}{c} \text{The innovations process is a sequence} \\[4pt] \text{of i.i.d. mean 0 and variance 1 r.v.'s.} \end{array}$$

Throughout this paper we assume that (1.1)–(1.4) hold, so that, by definition, $\mu_1 = 0$ and $\mu_2 = 1$.

The existence of a unique strictly stationary solution of (1.1) and (1.2) is well established. See [6, 7] for details. In this paper a minimal set of conditions in [6, 7] for the existence and stationarity of the GARCH$(p, q)$ sequence $\{X_t, -\infty < t < \infty\}$ is assumed, plus the assumption (1.4).

Estimation of the parameter $\boldsymbol{\theta} = (\alpha_0, \alpha_1, \ldots, \alpha_p, \beta_1, \ldots, \beta_q)$ has been investigated by several authors. We only cite those relevant to our investigation. Throughout this paper we assume that $\hat{\boldsymbol{\theta}}_n$ is an estimator of $\boldsymbol{\theta}$ based on a sample $X_0, X_1, \ldots, X_n$, and that it is $\sqrt{n}$ consistent, that is,

$$(1.5) \qquad\qquad \sqrt{n}|\hat{\boldsymbol{\theta}}_n - \boldsymbol{\theta}| = O_P(1),$$

where we use $|\cdot|$ to denote the maximum norm of vectors or matrices. Recently Berkes, Horváth and Kokoszka [3] studied the asymptotic properties of the quasi-maximum likelihood estimator for $\boldsymbol{\theta}$ in GARCH$(p, q)$ models under mild conditions. Berkes and Horváth [1] have shown that the quasi-maximum likelihood estimator cannot be $\sqrt{n}$ consistent if $E|\varepsilon_0|^k = \infty$ for some $0 < k < 4$. Hall and Yao [14] also studied inference for GARCH models under slightly stronger assumptions on the parameters. To remove such limitations as the finite fourth innovations moment, Berkes and Horváth [1] have used an arbitrary density function to replace the normal function used in the quasi-maximum likelihood and obtain the asymptotic properties (consistency and normality) under mild conditions. In particular, they show that the quasi-maximum likelihood estimator based on the standard symmetric exponential density function is $\sqrt{n}$ consistent only if $E\varepsilon_0^2 < \infty$.

The main goal of this paper is to construct high moment partial sum processes of residuals in a GARCH model. Since GARCH processes are defined in terms of their conditional variances, it is natural to construct model diagnostics in terms of sums or sums of squares of either the observed raw data or the driving noise estimates (residuals). As is well known from regression and linear time series, diagnostics based on residuals are often better tools. Sample skewness and kurtosis are based on sums of third and fourth powers, respectively, and thus, are also interesting diagnostic tools.

Usually the conditional variance $\sigma_t^2$ of (1.2) is estimated by

$$\hat{\sigma}_t^2 = \hat{\alpha}_0 + \sum_{i=1}^p \hat{\alpha}_i X_{t-i}^2 + \sum_{j=1}^q \hat{\beta}_j \hat{\sigma}_{t-j}^2, \qquad R = \max(p, q) \le t \le n.$$

One problem with the above estimation is that the initial conditional variance estimates $\hat{\sigma}_{R-1}^2, \hat{\sigma}_{R-2}^2, \ldots, \hat{\sigma}_{R-q}^2$ must be given. Hall and Yao [14] give an infinite-order moving average representation for $\sigma_t^2$ under the condition that $\sum_{i=1}^p \alpha_i + \sum_{j=1}^q \beta_j < 1$. A recursive representation of $\sigma_t^2$ is obtained by Berkes, Horváth and Kokoszka [3] under the weaker condition (than [14]) that $\sum_{j=1}^q \beta_j < 1$, which we are already assuming for the stationarity of $\{X_t, -\infty < t < \infty\}$ under the conditions of Bougerol and Picard [6, 7]. We use these later results of Berkes, Horváth and Kokoszka [3] to construct $\hat{\sigma}_t^2$, adapt their notation and conditions and give them here in some detail.

Let $\mathbf{u} = (\mathbf{s}, \mathbf{t}) \in \mathbb{R}^{p+q+1}$, $\mathbf{s} = (s_0, s_1, \ldots, s_p) \in \mathbb{R}^{p+1}$ and $\mathbf{t} = (t_1, \ldots, t_q) \in \mathbb{R}^q$. Define $c_i(\mathbf{u})$, $i = 1, 2, \ldots, R = \max\{q, p\}$ by: if $q \ge p$, then

$$c_0(\mathbf{u}) = s_0/(1 - (t_1 + \cdots + t_q)),$$

$$c_1(\mathbf{u}) = s_1,$$

$$c_2(\mathbf{u}) = s_2 + t_1 c_1(\mathbf{u}),$$

$$\vdots$$

$$c_p(\mathbf{u}) = s_p + t_1 c_{p-1}(\mathbf{u}) + \cdots + t_{p-1} c_1(\mathbf{u}),$$

$$c_{p+1}(\mathbf{u}) = t_1 c_p(\mathbf{u}) + \cdots + t_p c_1(\mathbf{u}),$$

$$\vdots$$

$$c_q(\mathbf{u}) = t_1 c_{q-1}(\mathbf{u}) + \cdots + t_{q-1} c_1(\mathbf{u}),$$

and if $q < p$, the equations above are replaced with

$$c_0(\mathbf{u}) = s_0/(1 - (t_1 + \cdots + t_q)),$$

$$c_1(\mathbf{u}) = s_1,$$

$$c_2(\mathbf{u}) = s_2 + t_1 c_1(\mathbf{u}),$$

$$\vdots$$

$$c_{q+1}(\mathbf{u}) = s_{q+1} + t_1 c_q(\mathbf{u}) + \cdots + t_q c_1(\mathbf{u}),$$

$$\vdots$$

$$c_p(\mathbf{u}) = s_p + t_1 c_{p-1}(\mathbf{u}) + \cdots + t_q c_{p-q}(\mathbf{u}).$$

For $i > R$, define

$$c_i(\mathbf{u}) = t_1 c_{i-1}(\mathbf{u}) + t_2 c_{i-2}(\mathbf{u}) + \cdots + t_q c_{i-q}(\mathbf{u}).$$

Let $0 < \underline{u} < \bar{u}$, $0 < \rho_0 < 1$, $q\underline{u} < \rho_0$, and define the parameter space as

$$\boldsymbol{\Theta} = \{\mathbf{u} : t_1 + \cdots + t_q \leq \rho_0, \underline{u} \leq \min(\mathbf{s}, \mathbf{t}) \leq \max(\mathbf{s}, \mathbf{t}) \leq \bar{u}\}.$$

In the rest of this paper we replace (1.3) with the stronger condition

(1.6)                    $\boldsymbol{\theta}$ is in the interior of $\boldsymbol{\Theta}$.

Now we are ready to give the recursive representation of the conditional variances by previous observations as given by Berkes, Horváth and Kokoszka [3]. Define

$$\sigma_t^2(\mathbf{u}) = c_0(\mathbf{u}) + \sum_{i=1}^{\infty} c_i(\mathbf{u}) X_{t-i}^2.$$

Then $\sigma_t^2(\mathbf{u})$ exists with probability one for all $\mathbf{u} \in \boldsymbol{\Theta}$. Also, $\sigma_t^2$ in (1.2) has the representation

(1.7)            $$\sigma_t^2 = \sigma_t^2(\boldsymbol{\theta}) = c_0(\boldsymbol{\theta}) + \sum_{i=1}^{\infty} c_i(\boldsymbol{\theta}) X_{t-i}^2.$$

Given $\hat{\boldsymbol{\theta}}_n$, we can estimate $\sigma_t^2(\boldsymbol{\theta})$ by

(1.8)            $$\tilde{\sigma}_t^2 = \sigma_t^2(\hat{\boldsymbol{\theta}}_n) = c_0(\hat{\boldsymbol{\theta}}_n) + \sum_{i=1}^{\infty} c_i(\hat{\boldsymbol{\theta}}_n) X_{t-i}^2.$$

In practice, we observe only $X_0, X_1, \ldots, X_n$. Hence, we use a truncated form and define, for $1 \leq t \leq n$,

(1.9)            $$\hat{\sigma}_t^2 = c_0(\hat{\boldsymbol{\theta}}_n) + \sum_{i=1}^{t} c_i(\hat{\boldsymbol{\theta}}_n) X_{t-i}^2.$$

Thus, the residual at time $t$ is

(1.10)            $$\hat{\varepsilon}_t = \frac{X_t}{\hat{\sigma}_t}, \qquad 1 \leq t \leq n.$$

The $k$th ($k = 1, 2, 3, 4, \ldots$) order moment partial sum process of residuals is defined as

$$(1.11) \qquad \hat{S}_n^{(k)}(u) = \sum_{t=1}^{[nu]} \hat{\varepsilon}_t^k, \qquad 0 \leq u \leq 1.$$

Its counterpart based on the i.i.d. innovations is defined as

$$(1.12) \qquad S_n^{(k)}(u) = \sum_{t=1}^{[nu]} \varepsilon_t^k, \qquad 0 \leq u \leq 1.$$

Collectively in $k$, we refer to these as high moment partial sum processes.

Denote

$$\partial c_i(\mathbf{u}) = \left( \frac{\partial c_i(\mathbf{u})}{\partial \alpha_0}, \frac{\partial c_i(\mathbf{u})}{\partial \alpha_1}, \ldots, \frac{\partial c_i(\mathbf{u})}{\partial \alpha_p}, \frac{\partial c_i(\mathbf{u})}{\partial \beta_1}, \ldots, \frac{\partial c_i(\mathbf{u})}{\partial \beta_q} \right) \in \mathbb{R}^{p+q+1},$$

$$\partial \log \sigma_t^2(\mathbf{u}) = \frac{\partial \sigma_t^2(\mathbf{u})}{\sigma_t^2(\mathbf{u})} = \frac{\partial c_0(\mathbf{u}) + \sum_{i=1}^{\infty} \partial c_i(\mathbf{u}) X_{t-i}^2}{c_0(\mathbf{u}) + \sum_{i=1}^{\infty} c_i(\mathbf{u}) X_{t-i}^2} \in \mathbb{R}^{p+q+1}$$

and

$$\psi(\boldsymbol{\theta}) = E(\partial \log \sigma_0^2(\boldsymbol{\theta})),$$

where $\partial(\cdot)$ is used as a shorthand for $\partial(\cdot)/\partial \mathbf{u}$ for convenience of writing.

We need two more regularity conditions in order to state our results:

$$(1.13) \qquad \varepsilon_0 \text{ is a nondegerate random variable}$$

and

$$(1.14) \qquad \lim_{x \to 0} x^{-\zeta} P\{|\varepsilon_0| \leq x\} = 0 \qquad \text{for some } \zeta > 0.$$

Lemma 3.1, (1.6), (1.14) and $E|\varepsilon_0|^\delta < \infty$ for some $\delta > 0$ imply the existence of $\psi(\boldsymbol{\theta})$.

THEOREM 1.1. *Suppose* (1.5), (1.6) *and* (1.14) *hold, and let* $k \geq 1$ *be an integer. If* $E|\varepsilon_0|^k < \infty$, *then*

$$\sup_{0 \leq u \leq 1} \left| \frac{1}{\sqrt{n}} (\hat{S}_n^{(k)}(u) - S_n^{(k)}(u)) + \frac{ku\mu_k}{2} \langle \psi(\boldsymbol{\theta}), \sqrt{n}(\hat{\boldsymbol{\theta}}_n - \boldsymbol{\theta}) \rangle \right| = o_P(1),$$

*where* $\langle \mathbf{x}, \mathbf{y} \rangle$ *is the inner product of the vectors* $\mathbf{x}$ *and* $\mathbf{y}$.

REMARK 1.1. Theorem 1.1 shows that the asymptotic properties of the high moment partial sum process $\{\hat{S}_n^{(k)}(u), 0 \leq u \leq 1\}$ depend on the parameters of the model unless $\mu_k = 0$, which can only happen if $k$ is an odd integer, i.e., not if $k$ is even. Recall that $\mu_1 = 0$ by (1.4). Thus, the ordinary partial sum process $\{\hat{S}_n^{(1)}(u), 0 \leq u \leq 1\}$ behaves as though the residuals $\{\hat{\varepsilon}_t, 1 \leq t \leq n\}$ were asymptotically the same as the unobservable innovations $\{\varepsilon_t, 1 \leq t \leq n\}$.

By Theorem 1.1, we immediately obtain the following CUSUM result, Theorem 1.2. It implies that the CUSUM normalized high moment partial sum process $\{\hat{S}_n^{(k)}(u) - u\hat{S}_n^{(k)}(1), 0 \le u \le 1\}$ behaves as though the residuals $\{\hat{\varepsilon}_t, 1 \le t \le n\}$ were asymptotically the same as the unobservable innovations $\{\varepsilon_t, 1 \le t \le n\}$.

THEOREM 1.2. *Let $k \ge 1$ be an integer and suppose that* (1.5), (1.6) *and* (1.14) *hold. If $E|\varepsilon_0|^k < \infty$, then*

$$\sup_{0 \le u \le 1} \frac{1}{\sqrt{n}}|(\hat{S}_n^{(k)}(u) - u\hat{S}_n^{(k)}(1)) - (S_n^{(k)}(u) - uS_n^{(k)}(1))| = o_P(1).$$

Let $\nu_k^2 = E(\varepsilon_0^k - \mu_k)^2 < \infty$. Then the invariance principle for partial sums for an i.i.d. sequence $\{\varepsilon_t^k\}$ (see, e.g., [4]) implies that

$$\left\{ \frac{S_n^{(k)}(u) - uS_n^{(k)}(1)}{\nu_k\sqrt{n}}, 0 \le u \le 1 \right\}$$

converges weakly in the Skorokhod space $D[0, 1]$ with $J_1$ topology to a Brownian bridge $\{B_0(u), 0 \le u \le 1\}$. Note that any topology of weak convergence that yields the invariance principle above could have been used here, but for definiteness we state that the $J_1$ topology is used.

The next result follows immediately from either Theorem 1.1 or Theorem 1.2.

COROLLARY 1.1. *Suppose* (1.5), (1.6), (1.13) *and* (1.14) *hold. If $E|\varepsilon_0|^{2k} < \infty$ for some integer $k \ge 1$, then*

$$\left\{ \frac{\hat{S}_n^{(k)}(u) - u\hat{S}_n^{(k)}(1)}{\nu_k\sqrt{n}}, \ 0 \le u \le 1 \right\}$$

*converges weakly in the Skorokhod space $D[0, 1]$ with $J_1$ topology to a Brownian bridge $\{B_0(u), 0 \le u \le 1\}$.*

REMARK 1.2. To use Corollary 1.1 for CUSUM tests of structural change of GARCH models, one needs to estimate $\nu_k$. The details are left to the next section.

THEOREM 1.3. *Suppose* (1.5), (1.6) *and* (1.14) *hold. If $E|\varepsilon_0|^k < \infty$ for an integer $k \ge 1$, then*

$$\left| \frac{1}{\sqrt{n}} \sum_{t=1}^n |\hat{\varepsilon}_t^k - \varepsilon_t^k| - \frac{k}{2}\psi_k(\sqrt{n}(\hat{\boldsymbol{\theta}}_n - \boldsymbol{\theta})) \right| = o_P(1),$$

*where*

$$\psi_k(\mathbf{u}) = E|\langle \varepsilon_0^k \, \partial \log \sigma_0^2(\boldsymbol{\theta}), \mathbf{u} \rangle|, \qquad \mathbf{u} \in \mathbb{R}^{p+q+1}.$$

REMARK 1.3.    The above theorem shows that the sum of the absolute deviations $|\hat{\varepsilon}_t^k - \varepsilon_t^k|$ depends on the parameters of the model. The existence of $\psi_k(\mathbf{u})$ follows from

$$E|\langle \varepsilon_0^k \, \partial \log \sigma_0^2(\boldsymbol{\theta}), \mathbf{u} \rangle| \leq |\mathbf{u}| E|\varepsilon_0|^k E|\partial \log \sigma_0^2(\boldsymbol{\theta})| < \infty$$

since $\varepsilon_0$ and $\partial \log \sigma_0^2(\boldsymbol{\theta})$ are independent and $E|\partial \log \sigma_0^2(\boldsymbol{\theta})| < \infty$ by Lemma 3.1. It is easy to check that $\psi_k(\mathbf{u})$ is symmetric about 0 and is Lipschitz by

$$|\psi_k(\mathbf{u}) - \psi_k(\mathbf{u}^*)| \leq |\mathbf{u} - \mathbf{u}^*| E|\varepsilon_0|^k E|\partial \log \sigma_0^2(\boldsymbol{\theta})| \qquad \forall \mathbf{u}, \mathbf{u}^* \in \mathbb{R}^{p+q+1},$$

where $\mathbf{u}^* = (\mathbf{s}^*, \mathbf{t}^*)$, $\mathbf{s}^* = (s_0^*, s_1^*, \ldots, s_p^*) \in \mathbb{R}^{p+1}$ and $\mathbf{t}^* = (t_1^*, \ldots, t_q^*) \in \mathbb{R}^q$.

Before formulating the next result, we need to modify the high moment partial sum processes of (1.11) and (1.12). The $k$th-order moment residual centered partial sum process is defined as

$$(1.15) \qquad \hat{T}_n^{(k)}(u) = \sum_{t=1}^{[nu]} (\hat{\varepsilon}_t - \bar{\hat{\varepsilon}})^k, \qquad 0 \leq u \leq 1,$$

where $\bar{\hat{\varepsilon}}$ is the sample mean of the residuals. Its counterpart based on the i.i.d. innovations is

$$T_n^{(k)}(u) = \sum_{t=1}^{[nu]} (\varepsilon_t - \bar{\varepsilon})^k, \qquad 0 \leq u \leq 1,$$

where $\bar{\varepsilon}$ is the sample mean of innovations.

Obviously, $\hat{\sigma}_{(n)}^2 = \hat{T}_n^{(2)}(1)/n$ is the sample moment estimator of $\mu_2$. In fact, by (3.10) $\hat{\sigma}_{(n)}^2 \to \mu_2$ in probability under the minimal condition $\mu_2 < \infty$. Since $\mu_2 = 1$, then $\hat{\sigma}_{(n)}^2$ may seem to be an unnecessary estimator. However, it will play an important role when it is used to self-normalize $\hat{T}_n^{(k)}(u)$. Denote $\sigma_{(n)}^2 = T_n^{(2)}(1)/n$, and note that it is the sample variance of the true innovations, except with divisor $n$ instead of $n-1$.

THEOREM 1.4.    *Suppose* (1.5), (1.6), (1.13) *and* (1.14) *hold and* $k \geq 1$ *is an integer. If* $E|\varepsilon_0|^{\max\{k,2\}} < \infty$, *then*

$$\sup_{0 \leq u \leq 1} \frac{1}{\sqrt{n}} \left| \frac{\hat{T}_n^{(k)}(u)}{\hat{\sigma}_{(n)}^k} - \frac{T_n^{(k)}(u)}{\sigma_{(n)}^k} \right| = o_P(1).$$

REMARK 1.4.    Theorem 1.4 implies that the self-normalized (or estimated scale normalized) high moment centered partial sum process $\{\hat{T}_n^{(k)}(u)/\hat{\sigma}_{(n)}^k, 0 \leq u \leq 1\}$ behaves as though the residuals $\{\hat{\varepsilon}_t, 1 \leq t \leq n\}$ were asymptotically the same as the unobservable innovations $\{\varepsilon_t, 1 \leq t \leq n\}$.

REMARK 1.5.    When $k = 1$, it is easy to verify that

$$\sup_{0 \leq u \leq 1} |\hat{T}_n^{(1)}(u) - (\hat{S}_n^{(1)}(u) - u\hat{S}_n^{(1)}(1))| \leq |\bar{\hat{\varepsilon}}|.$$

Hence, given $E\varepsilon_0^2 < \infty$, Theorem 1.1 (cf. Remark 1.1) implies $\bar{\hat{\varepsilon}} = O_P(1/\sqrt{n})$ and Corollary 1.1 implies that

$$\left\{ \frac{\hat{T}_n^{(1)}(u)}{\hat{\sigma}_{(n)}\sqrt{n}}, 0 \leq u \leq 1 \right\}$$

converges weakly in the Skorokhod space $D[0,1]$ to a Brownian bridge $\{B_0(u), 0 \leq u \leq 1\}$.

Let $\lambda_k = \mu_k/\mu_2^{k/2}$ for $k \geq 1$ and define $\lambda_0 = 1$. For each $k \geq 1$, let $\{B^{(k)}(u), 0 \leq u \leq 1\}$ be a zero mean Gaussian process with covariance

$$\begin{aligned}
EB^{(k)}(u)B^{(k)}(v) &= (\lambda_{2k} - \lambda_k^2)(u \wedge v) \\
&\quad + k\lambda_{k-1}(k\lambda_{k-1} + k\lambda_k\lambda_3 - 2\lambda_{k+1})uv \\
&\quad + k\lambda_k((1 - k/4)\lambda_k + k\lambda_k\lambda_4/4 - \lambda_{k+2})uv
\end{aligned} \tag{1.16}$$

for any $0 \leq u, v \leq 1$, where $u \wedge v = \min(u, v)$.

If $\mu_{2k} < \infty$, then Lemma 3.8 implies

$$\left\{ \frac{1}{\sqrt{n}} \left( \frac{T_n^{(k)}(u)}{\sigma_{(n)}^k} - nu\lambda_k \right), 0 \leq u \leq 1 \right\}$$

converges weakly to the Gaussian process $\{B^{(k)}(u), 0 \leq u \leq 1\}$. By Theorem 1.4, we immediately obtain the following corollary.

COROLLARY 1.2.    *If (1.5), (1.6), (1.13) and (1.14) hold, then $E|\varepsilon_0|^{2k} < \infty$ for some integer $k \geq 1$ implies that*

$$\left\{ \frac{1}{\sqrt{n}} \left( \frac{\hat{T}_n^{(k)}(u)}{\hat{\sigma}_{(n)}^k} - nu\lambda_k \right), 0 \leq u \leq 1 \right\}$$

*converges weakly to the Gaussian process $\{B^{(k)}(u), 0 \leq u \leq 1\}$.*

REMARK 1.6.    If the innovation distribution is symmetric about 0, then the covariance in (1.16) can be simplified. If $k$ is odd, then

$$EB^{(k)}(u)B^{(k)}(v) = \lambda_{2k}(u \wedge v) + k\lambda_{k-1}(k\lambda_{k-1} - 2\lambda_{k+1})uv$$

and if $k$ is even, then

$$EB^{(k)}(u)B^{(k)}(v) = (\lambda_{2k} - \lambda_k^2)(u \wedge v) + k\lambda_k((1 - k/4)\lambda_k + k\lambda_k\lambda_4/4 - \lambda_{k+2})uv.$$

REMARK 1.7. Based on the facts that $\lambda_0 = 1$, $\lambda_1 = 0$ and $\lambda_2 = 1$, (1.16) becomes $EB^{(1)}(u)B^{(1)}(v) = u \wedge v - uv$ for any $0 \le u, v \le 1$. That is, $\{B^{(1)}(u), 0 \le u \le 1\}$ is a Brownian bridge. Hence, the result of Corollary 1.2 for $k = 1$ matches that in Remark 1.5. For $k = 2$, simple calculations from (1.16) show that $EB^{(2)}(u)B^{(2)}(v) = (\lambda_4 - 1)(u \wedge v - uv)$ for any $0 \le u, v \le 1$. Notice that $\nu_2 = \mu_4 - \mu_2^2 = \mu_2^2(\lambda_4 - 1)$. Thus, when $k = 2$, the result of Corollary 1.2 matches that of Corollary 1.1. In general, the Gaussian process $\{B^{(k)}(u), 0 \le u \le 1\}$ for $k \ge 3$ depends on the moments of the innovation distribution and cannot be identified with a specific classic process such as a Brownian motion or Brownian bridge.

Notice that Corollary 1.2 gives the weak convergence of the self-normalized high moment centered partial sum process $\{\hat{T}_n^{(k)}(u)/\hat{\sigma}_{(n)}^k, 0 \le u \le 1\}$ of residuals for a fixed $k$. The following result considers the joint weak convergence of two self-normalized high moment centered partial sum processes.

THEOREM 1.5. Assume that (1.5), (1.6), (1.13) and (1.14) hold. Assume also that $k \ge 1$ is an odd number and $\mu_3 = \mu_k = \mu_{k+2} = \mu_{2k+1} = 0$. Then $E|\varepsilon_0|^{2(k+1)} < \infty$ implies that

$$\left\{ \frac{1}{\sqrt{n}} \left( \frac{\hat{T}_n^{(k)}(u)}{\hat{\sigma}_{(n)}^k} - nu\lambda_k, \frac{\hat{T}_n^{(k+1)}(v)}{\hat{\sigma}_{(n)}^{k+1}} - nu\lambda_{k+1} \right), 0 \le u, v \le 1 \right\}$$

converges weakly, in the Skorokhod space $D^2[0,1]$ equipped with the product $J_1$ topology, to a two-dimensional Gaussian process $\{(B^{(k)}(u), B^{(k+1)}(v)), 0 \le u, v \le 1\}$, where $\{B^{(k)}(u), 0 \le u \le 1\}$ and $\{B^{(k+1)}(v), 0 \le v \le 1\}$ are two independent Gaussian processes defined by (1.16).

REMARK 1.8. The conditions $\mu_3 = \mu_k = \mu_{k+2} = \mu_{2k+1} = 0$ can be replaced with the stronger condition that the innovation distribution is symmetric about 0.

## 2. Applications.

This section considers applications of the high moment residual partial sums to a change-point problem, goodness-of-fit and the construction of a kernel density estimate of the unobservable innovation distribution.

### 2.1. CUSUM tests for structural change of GARCH models.

In this subsection we consider the CUSUM normalized high moment partial sum process $\{\hat{S}_n^{(k)}(u) - u\hat{S}_n^{(k)}(1), 0 \le u \le 1\}$ defined in Theorem 1.2. It is related to the standard CUSUM test introduced by Brown, Durbin and Evans [9], which was one of the first tests on structural change with unknown break point.

We first consider a structural change in the conditional mean for GARCH models. We can formulate it as the following hypothesis test. The null hypothesis is "no-change in the conditional mean,"

$$(2.1) \qquad H_0 : \begin{cases} X_t = \sigma_t \varepsilon_t, \\ \sigma_t^2 = \alpha_0 + \sum_{i=1}^{p} \alpha_i X_{t-i}^2 + \sum_{j=1}^{q} \beta_j \end{cases}, \qquad t = 0, 1, \ldots, n,$$

and the alternative is "one change in the conditional mean,"

$$H_a : \begin{cases} X_t = \sigma_t \varepsilon_t, \\ \sigma_t^2 = \alpha_0 + \sum_{i=1}^{p} \alpha_i X_{t-i}^2 + \sum_{j=1}^{q} \beta_j \sigma_{t-j}^2 \end{cases}, \qquad t = 0, \ldots, [nu^*], \\ \begin{cases} X_t = \sigma_t \varepsilon_t + \mu, \\ \sigma_t^2 = \alpha_0 + \sum_{i=1}^{p} \alpha_i (X_{t-i} - \mu)^2 + \sum_{j=1}^{q} \beta_j \sigma_{t-j}^2 \end{cases}, \qquad t = [nu^*]+1, \ldots, n,$$

where $\mu \neq 0$ and $0 < u^* < 1$.

To test the above hypothesis, we use the standard CUSUM test constructed from residuals as

$$CUSUM^{(1)} = \max_{1 \leq i < n} \frac{|\sum_{t=1}^{i} \hat{\varepsilon}_t - i\hat{\bar{\varepsilon}}|}{\hat{\sigma}_{(n)} \sqrt{n}}.$$

By a straightforward calculation, it is easy to verify that

$$CUSUM^{(1)} = \sup_{0 \leq u \leq 1} \frac{|\hat{S}_n^{(1)}(u) - u\hat{S}_n^{(1)}(1)|}{\hat{\sigma}_{(n)} \sqrt{n}} + o_P(1),$$

provided that $E\varepsilon_0^2 < \infty$. Therefore, by Corollary 1.1, under $H_0$,

$$CUSUM^{(1)} \xrightarrow{\mathcal{D}} \sup_{0 \leq u \leq 1} |B_0(u)|,$$

where $\{B_0(u), 0 \leq u \leq 1\}$ is a Brownian bridge. Hence, we can reject the $H_0$ in favor of $H_a$ if $CUSUM^{(1)}$ is large.

REMARK 2.1. The statistic $CUSUM^{(1)}$ involves the estimator $\sqrt{\hat{\mu}_2}$ of $\sqrt{\mu_2}$. However, according to the GARCH model setup, $\mu_2 = 1$. Thus, the term $\hat{\sigma}_{(n)}$ can be dropped in $CUSUM^{(1)}$ to obtain a related test statistic.

A second interesting structural change hypothesis concerns a change in the conditional variance of a GARCH model. We use the above $H_0$ as the null hypothesis for "no-change in the conditional variance" against the "one change in the conditional variance" alternative

$$H_{a'}: \begin{cases} X_t = \sigma_t \varepsilon_t, \\ \sigma_t^2 = \begin{cases} \alpha_0 + \sum_{i=1}^p \alpha_i X_{t-i}^2 + \sum_{j=1}^q \beta_j \sigma_{t-j}^2, & \text{if } t = 0, \dots, [nu^*], \\ \alpha_0' + \sum_{i=1}^p \alpha_i' X_{t-i}^2 + \sum_{j=1}^q \beta_j' \sigma_{t-j}^2, & \text{if } t = [nu^*]+1, \dots, n, \end{cases} \end{cases}$$

where

$$(\alpha_0, \alpha_1, \dots, \alpha_p, \beta_1, \dots, \beta_q) \neq (\alpha_0', \alpha_1', \dots, \alpha_p', \beta_1', \dots, \beta_q')$$

and

$$0 < u^* < 1.$$

In the following we propose two CUSUM statistics. The first statistic is defined as

$$CUSUM_1^{(2)} = \max_{1 \leq i < n} \frac{|\sum_{t=1}^i \hat{\varepsilon}_t^2 - i \sum_{t=1}^n \hat{\varepsilon}_t^2 / n|}{\hat{\nu}_2 \sqrt{n}},$$

where

$$\hat{\nu}_2^2 = \frac{1}{n} \sum_{t=1}^n ((\hat{\varepsilon}_t - \bar{\hat{\varepsilon}})^2 - \hat{\sigma}_{(n)}^2)^2$$

is an estimator of $\nu_2 = E(\varepsilon_0^2 - \mu_2)^2$. The statistic $\hat{\nu}_2$ uses the fact that $\lambda_2 = 1$ is known from the definition of the GARCH process and is not estimated; see (1.4). The second statistic is defined as

$$CUSUM_2^{(2)} = \max_{1 \leq i < n} \frac{|\sum_{t=1}^i (\hat{\varepsilon}_t - \bar{\hat{\varepsilon}})^2 - i \hat{\sigma}_{(n)}^2|}{\hat{\nu}_2 \sqrt{n}},$$

that is, $CUSUM_2^{(2)}$ is centered about the residual sample mean $\bar{\hat{\varepsilon}}$ in contrast to the no centering $CUSUM_1^{(2)}$. Again, by straightforward calculations, it is easy to show that

$$CUSUM_1^{(2)} = \sup_{0 \leq u \leq 1} \frac{|\hat{S}_n^{(2)}(u) - u \hat{S}_n^{(2)}(1)|}{\hat{\nu}_2 \sqrt{n}} + o_P(1)$$

and

$$CUSUM_2^{(2)} = \sup_{0 \leq u \leq 1} \frac{\hat{\sigma}_{(n)}^2}{\hat{\nu}_2 \sqrt{n}} \left| \frac{\hat{T}_n^{(2)}(u)}{\hat{\sigma}_{(n)}^2} - nu\lambda_2 \right| + o_P(1),$$

TABLE 1
*Size and power of $CUSUM_2^{(2)}$ statistic for GARCH(1, 1)*

| $\varepsilon_0 = N(0, 1)$ | $n = 500$ | $n = 1000$ | $n = 1500$ | $n = 3000$ |
|---|---|---|---|---|
| $\boldsymbol{\theta}$ = (0.0002, 0.1, 0.7) (Null) | 0.0230 | 0.0358 | 0.0394 | 0.0412 |
| $\boldsymbol{\theta}'$ = (0.0003, 0.1, 0.7) | 0.2342 | 0.6484 | 0.8752 | 0.9958 |
| $\boldsymbol{\theta}'$ = (0.0002, 0.167, 0.7) | 0.1316 | 0.3908 | 0.5998 | 0.9186 |
| $\boldsymbol{\theta}'$ = (0.0002, 0.1, 0.767) | 0.1792 | 0.5470 | 0.8264 | 0.9924 |
| $\boldsymbol{\theta}$ = (0.0002, 0.1, 0.8) (Null) | 0.0162 | 0.0320 | 0.0370 | 0.0386 |
| $\boldsymbol{\theta}'$ = (0.0003, 0.1, 0.8) | 0.1040 | 0.3922 | 0.6570 | 0.9642 |
| $\boldsymbol{\theta}'$ = (0.0002, 0.167, 0.8) | 0.1840 | 0.6260 | 0.8786 | 0.9978 |
| $\boldsymbol{\theta}'$ = (0.0002, 0.1, 0.867) | 0.1768 | 0.5980 | 0.9320 | 1.0000 |

provided that $E\varepsilon_0^4 < \infty$. Therefore, by Corollaries 1.1 and 1.2 (cf. Remark 1.7), under $H_0$,

$$CUSUM_i^{(2)} \xrightarrow{\mathcal{D}} \sup_{0 \leq u \leq 1} |B_0(u)|, \qquad i = 1, 2,$$

where $\{B_0(u), 0 \leq u \leq 1\}$ is a Brownian bridge. Hence, we can reject the $H_0$ in favor of $H_{a'}$ whenever $CUSUM_i^{(2)}$ ($i = 1, 2$) is large.

Preliminary empirical studies show promising results from the above proposed test statistics. They outperform the CUSUM test constructed from the squares of the original data by Kim, Cho and Lee [17]. Some details are given in [23]. Independently, Kokoszka and Leipus [18] also study a change point for an ARCH process, again based on the original observations and not residuals. Here we list empirical sizes and powers of the $CUSUM_2^{(2)}$ test. The significance level is 5% with 1.358 as the critical value, the break point $u^*$ at the alternative is 0.5, and the number of replicates is 5000. Tables 1 and 2 show that there are size distortions of the $CUSUM_2^{(2)}$ test, but less serious with large sample size. The null and alternatives given in these tables are slightly different from the cases studied in [17]. Their tables use $\alpha_0$ as 0.2 or 0.3. When fitting GARCH models to financial stock returns data, typically a much smaller value of $\alpha_0$ is found and, hence, our tables use values of $\alpha_0$ such as 0.0002 and 0.0003. A simulation with $\alpha_0$ as 0.2 or 0.3 was also undertaken, but not reported here. In addition, we include the near integrated GARCH cases with $\alpha_1 = 0.1$ or 0.167 and $\beta_1 = 0.8$ or 0.867 that were shown to have poor performance in Kim, Cho and Lee's [17]. The $CUSUM_2^{(2)}$ test outperforms Kim, Cho and Lee's test in all instances, with both large and small values of $\alpha_0$ and, in particular, it has substantial power gains when the innovation distribution is $t(8)$.

2.2. *Jarque–Bera normality test.* In this subsection we consider the self-normalized high moment centered partial sum process (1.15) for $k = 3$ and $k = 4$. They correspond to the sample skewness partial sum process as

$$\hat{\gamma}_n(u) = \frac{\hat{T}_n^{(3)}(u)/n}{\hat{\sigma}_{(n)}^3}, \qquad 0 \leq u \leq 1,$$

and the sample kurtosis process as

$$\hat{\kappa}_n(u) = \frac{\hat{T}_n^{(4)}(u)/n}{\hat{\sigma}_{(n)}^4}, \qquad 0 \leq u \leq 1.$$

The sample skewness and kurtosis of the residuals are $\hat{\gamma}_n(1)$ and $\hat{\kappa}_n(1)$, respectively. Omnibus statistics based on sample skewness and kurtosis have been used to test normality. Bowman and Shenton [8] and Gasser [13] give details of this method. The basic idea is to construct the statistic

$$\frac{n}{\sigma_\gamma^2}(\hat{\gamma}_n(1) - \lambda_3)^2 + \frac{n}{\sigma_\kappa^2}(\hat{\kappa}_n(1) - \lambda_4)^2,$$

where, by (1.16),

$$\sigma_\gamma^2 = E(B^{(3)}(1))^2 = (\lambda_6 - \lambda_3^2) + 3(3 + 3\lambda_3^2 - 2\lambda_4) + 3\lambda_3(\lambda_3/4 + 3\lambda_3\lambda_4/4 - \lambda_5)$$

and

$$\sigma_\kappa^2 = E(B^{(4)}(1))^2 = (\lambda_8 - \lambda_4^2) + 4\lambda_3(4\lambda_3 + 4\lambda_3\lambda_4 - 2\lambda_5) + 4\lambda_4(\lambda_4^2 - \lambda_6).$$

Assume that the innovation distribution is symmetric about 0. Then by Theorem 1.5

$$(2.2) \qquad \frac{n}{\sigma_\gamma^2}(\hat{\gamma}_n(1) - \lambda_3)^2 + \frac{n}{\sigma_\kappa^2}(\hat{\kappa}_n(1) - \lambda_4)^2 \xrightarrow{\mathcal{D}} \chi_{(2)}^2.$$

TABLE 2
*Size and Power of $CUSUM_2^{(2)}$ statistic for GARCH(1, 1)*

| $\varepsilon_0 = t(8)$ | $n = 500$ | $n = 1000$ | $n = 1500$ | $n = 3000$ |
|---|---|---|---|---|
| $\boldsymbol{\theta} = (0.0002, 0.1, 0.7)$ (Null) | 0.0234 | 0.0302 | 0.0336 | 0.0396 |
| $\boldsymbol{\theta}' = (0.0003, 0.1, 0.7)$ | 0.1524 | 0.4286 | 0.6708 | 0.9540 |
| $\boldsymbol{\theta}' = (0.0002, 0.167, 0.7)$ | 0.0752 | 0.1986 | 0.3234 | 0.6620 |
| $\boldsymbol{\theta}' = (0.0002, 0.1, 0.767)$ | 0.1056 | 0.3370 | 0.5542 | 0.9132 |
| $\boldsymbol{\theta} = (0.0002, 0.1, 0.8)$ (Null) | 0.0188 | 0.0256 | 0.0318 | 0.0378 |
| $\boldsymbol{\theta}' = (0.0003, 0.1, 0.8)$ | 0.0716 | 0.2432 | 0.4354 | 0.8126 |
| $\boldsymbol{\theta}' = (0.0002, 0.167, 0.8)$ | 0.0932 | 0.3450 | 0.5758 | 0.9308 |
| $\boldsymbol{\theta}' = (0.0002, 0.1, 0.867)$ | 0.1126 | 0.3770 | 0.7198 | 0.9924 |

In the special case where the innovation distribution is standard normal, for which $\lambda_3 = 0$, $\lambda_4 = 3$, $\sigma_\gamma^2 = 6$ and $\sigma_\kappa^2 = 24$, then (2.2) becomes the Jarque–Bera (JB) statistic and has a $\chi_{(2)}^2$ limit in distribution,

$$(2.3) \qquad JB = \frac{n}{6}\hat{\gamma}_n^2(1) + \frac{n}{24}(\hat{\kappa}_n(1) - 3)^2 \xrightarrow{\mathcal{D}} \chi_{(2)}^2.$$

The statistic $JB$ in (2.3) is the Jarque–Bera normality test widely used in econometrics and implemented in standard statistical packages such as S-PLUS, and Jarque and Bera [15] show that $JB$ is a Lagrange multiplier test statistic of normality against alternatives within the Pearson family of distributions, which includes the beta, gamma and Student's $t$ distributions among others. They point out that it is asymptotically equivalent to the likelihood ratio test, implying it has the same asymptotic power characteristics and, hence, has maximum local asymptotic power [11]. Therefore, a test based on $JB$ is asymptotically locally most powerful against the Pearson family, and (2.3) shows that $JB$ is asymptotically distributed as $\chi^2(2)$. The hypothesis of normality is rejected for large sample size, if the computed value of $JB$ is greater than the appropriate critical value of a $\chi_{(2)}^2$. Lu [19] has used Monte Carlo simulation to obtain critical values for several different $n$. Based on this, a finite sample size correction can also be used to improve the choice of the critical value. Lu [19] obtained the finite sample size correction for the size 0.05 critical values of the $JB$ test as

$$JB_{0.05} = 5.991645 - 15.17n^{-1/2} + 345.9n^{-1} - 3110.8n^{-3/2}, \qquad n \geq 100.$$

We are unaware of any other results studying the Jarque–Bera test for GARCH residuals. Kilian and Demiroglu [16] studied the Jarque–Bera test for autoregressive residuals.

2.3. *Kernel density estimation of the innovation distribution.* The omnibus type statistic discussed in Section 2.2 provides a means to test a specific type of unobservable innovation distribution, such as normal, Student-$t$ and two-sided exponential. However, in practice, the normal innovation assumption is often rejected, as are other known types of distributions. Rather than focusing on identifying the innovation distribution to a specific member of a family, in this subsection we turn to a nonparametric kernel density estimation based on the residuals. This would be needed if one wished to implement a semi-parametric bootstrap methodology in this setting.

Assume that the innovation distribution has a uniformly continuous density function $f(x)$ which is unknown. Let $h_n$ be a sequence of positive numbers and $K$ be a probability density function (kernel) with mean 0 and variance 1. Then the kernel density estimation of $f(x)$ based on the residuals is defined as

$$\hat{f}_n(x) = \frac{1}{nh_n}\sum_{t=1}^{n} K\left(\frac{x - \hat{\varepsilon}_t}{h_n}\right), \qquad x \in \mathbb{R}.$$

Its counterpart based on i.i.d. innovations is defined as

$$f_n(x) = \frac{1}{nh_n} \sum_{t=1}^{n} K\left(\frac{x - \varepsilon_t}{h_n}\right), \qquad x \in \mathbb{R}.$$

THEOREM 2.1. *Assume that* (1.5), (1.6) *and* (1.14) *hold. In addition, we assume that the following three conditions hold:*

(i) $h_n > 0, h_n \to 0, \sqrt{n}h_n^2 \to \infty$;
(ii) $\sup_{|x|>b} |x|K(x) \to 0$ *as* $b \to \infty$;
(iii) $K$ *is Lipschitz, that is, there exists a constant* $C$ *such that*

$$|K(x) - K(y)| \le C|x - y| \qquad \forall\, x, y \in \mathbb{R}.$$

*Then* $E|\varepsilon_0| < \infty$ *implies that*

$$\sup_{x \in \mathbb{R}} |\hat{f}_n(x) - f_n(x)| = o_P(1).$$

The proof of Theorem 2.2 follows easily from Theorem 1.3. Given the conditions in Theorem 2.2, we have (cf. [21])

$$\sup_{x \in \mathbb{R}} |f_n(x) - f(x)| = o_P(1).$$

Thus, by Theorem 2.2,

$$\sup_{x \in \mathbb{R}} |\hat{f}_n(x) - f(x)| = o_P(1).$$

Notice in the above result that only the finite first innovation moment and $\sqrt{n}$ consistency of the parameter estimate are required.

**3. Proofs.** This section begins with a proof of Theorem 1.1. It is given in a sketch or overview form, with the details given in a series of lemmas which are placed in the later part of this section. The proofs of Theorems 1.3 and 1.4 rely on the proof of Theorem 1.1.

PROOF OF THEOREM 1.1. Let

$$\tilde{\varepsilon}_t = \frac{X_t}{\tilde{\sigma}_t}, \qquad 1 \le t \le n \quad \text{and} \quad \tilde{S}_n^{(k)}(u) = \sum_{t=1}^{[nu]} \tilde{\varepsilon}_t^k, \qquad 0 \le u \le 1,$$

where $\tilde{\sigma}_t$ is defined in (1.8). By (1.10) and (1.11),

$$\hat{\varepsilon}_t = \tilde{\varepsilon}_t\left(1 + \frac{\tilde{\sigma}_t - \hat{\sigma}_t}{\hat{\sigma}_t}\right)$$

and, hence,

$$\hat{S}_n^{(k)}(u) = \tilde{S}_n^{(k)}(u) + \sum_{i=1}^{k} \binom{k}{i} \sum_{t=1}^{[nu]} \tilde{\varepsilon}_t^k \left(\frac{\tilde{\sigma}_t - \hat{\sigma}_t}{\hat{\sigma}_t}\right)^i.$$

Thus, Theorem 1.1 follows if

$$(3.1) \quad \sup_{0 \le u \le 1} \left| \frac{1}{\sqrt{n}} (\tilde{S}_n^{(k)}(u) - S_n^{(k)}(u)) + \frac{ku\mu_k}{2} \langle \psi(\boldsymbol{\theta}), \sqrt{n}(\hat{\boldsymbol{\theta}}_n - \boldsymbol{\theta}) \rangle \right| = o_P(1)$$

and

$$(3.2) \quad \sum_{t=1}^{n} |\tilde{\varepsilon}_t|^k \left| \frac{\tilde{\sigma}_t - \hat{\sigma}_t}{\hat{\sigma}_t} \right|^i = O_P(1), \qquad 1 \le i \le k.$$

The sample conditional standard deviation estimates $\hat{\sigma}_t$ are uniformly bounded away from 0 in probability. This is argued as follows. Since $\boldsymbol{\theta}$ is $\sqrt{n}$-consistent, there exists an open ball in the interior of $\boldsymbol{\Theta}$ such that, for any small $\eta > 0$, then $\boldsymbol{\theta}$ belongs to this open ball with probability $\ge 1 - \eta$ as $n \to \infty$. Therefore, $\hat{\sigma}_t^2 \ge \hat{\alpha}_0 > \frac{1}{2}\alpha_0 > 0$ with probability $\ge 1 - \eta$. A stronger result than (3.2) is given in Lemma 3.5.

Let

$$(3.3) \quad g_t(\mathbf{u}) = \frac{\sqrt{n}(\sigma_t^2(\boldsymbol{\theta} + n^{-1/2}\mathbf{u}) - \sigma_t^2(\boldsymbol{\theta}))}{\sigma_t^2(\boldsymbol{\theta})}, \qquad \mathbf{u} \in \mathbb{R}^{p+q+1},$$

and

$$(3.4) \quad \varepsilon_t(\mathbf{u}) = \frac{\varepsilon_t}{\sqrt{1 + n^{-1/2}g_t(\mathbf{u})}}.$$

Though $g_t(\mathbf{u})$ depends on $n$, we omit it for convenience of notation.

Using (1.1), (1.7), (1.8), (3.3) and (3.4), we can rewrite $\tilde{\varepsilon}_t$ as

$$\tilde{\varepsilon}_t = \varepsilon_t(\sqrt{n}(\hat{\boldsymbol{\theta}}_n - \boldsymbol{\theta})), \qquad 1 \le t \le n,$$

that is, $\tilde{\varepsilon}_t = \varepsilon_t(\mathbf{u})$ with $\mathbf{u} = \sqrt{n}(\hat{\boldsymbol{\theta}}_n - \boldsymbol{\theta})$.

Hence, by (1.5), (1.11) and (1.12), we can prove (3.1) if, for any $b > 0$,

$$\sup_{|\mathbf{u}| \le b} \sup_{0 \le u \le 1} \left| \frac{1}{\sqrt{n}} \sum_{t=1}^{[nu]} (\varepsilon_t^k(\mathbf{u}) - \varepsilon_t^k) + \frac{ku\mu_k}{2} \langle \psi(\boldsymbol{\theta}), \mathbf{u} \rangle \right| = o_P(1).$$

This last part follows by

$$(3.5) \quad \sup_{|\mathbf{u}| \le b} \frac{1}{\sqrt{n}} \sum_{t=1}^{n} \left| \varepsilon_t^k(\mathbf{u}) - \varepsilon_t^k \left(1 - \frac{k}{2\sqrt{n}} \langle \partial \log \sigma_t^2(\boldsymbol{\theta}), \mathbf{u} \rangle \right) \right| = o_P(1)$$

and

$$(3.6) \quad \sup_{|\mathbf{u}| \le b} \sup_{0 \le u \le 1} \left| \frac{1}{n} \sum_{t=1}^{[nu]} \varepsilon_t^k \langle \partial \log \sigma_t^2(\boldsymbol{\theta}), \mathbf{u} \rangle - u\mu_k \langle \psi(\boldsymbol{\theta}), \mathbf{u} \rangle \right| = o_P(1).$$

The proof of (3.6) follows by taking $\sup_{|\mathbf{u}|\leq b}$ into the inner product first, then applying Lemma 3.6 and noting that $\langle\partial\log\sigma_t^2(\boldsymbol{\theta}),\mathbf{u}\rangle = h(\varepsilon_{t-1},\varepsilon_{t-2},\dots)$ for an appropriate function $h$. We also use the fact that $\varepsilon_t$ and $\partial\log\sigma_t^2(\boldsymbol{\theta})$ are independent, and that, by Lemma 3.1, $E(|\partial\log\sigma_t^2(\boldsymbol{\theta})|) < \infty$.

The main idea to prove (3.5) is to have a proper approximation of $1/\sqrt{1+n^{-1/2}g_t(\mathbf{u})}$ so that

$$\frac{1}{\sqrt{1+n^{-1/2}g_t(\mathbf{u})}} = 1 - \frac{g_t(\mathbf{u})}{2\sqrt{n}} + o_P\left(\frac{1}{\sqrt{n}}\right)$$

$$= 1 - \frac{\langle\partial\log\sigma_t^2(\boldsymbol{\theta}),\mathbf{u}\rangle}{2\sqrt{n}} + o_P\left(\frac{1}{\sqrt{n}}\right)$$

uniformly in $|\mathbf{u}|\leq b$ and $1\leq t\leq n$. Lemmas 3.1–3.4 are devoted to showing that this approximation holds.

We divide the proof of (3.5) into two parts. By (3.4), equation (3.5) will follow if

$$\sup_{|\mathbf{u}|\leq b}\frac{1}{\sqrt{n}}\sum_{t=1}^{n}|\varepsilon_t|^k\left|\left(\frac{1}{\sqrt{1+n^{-1/2}g_t(\mathbf{u})}}\right)^k - \left(1-\frac{\langle\partial\log\sigma_t^2(\boldsymbol{\theta}),\mathbf{u}\rangle}{2\sqrt{n}}\right)^k\right| = o_P(1)$$

and

$$(3.7)\quad\begin{aligned}\sup_{|\mathbf{u}|\leq b}\frac{1}{\sqrt{n}}\sum_{t=1}^{n}|\varepsilon_t|^k&\left|\left(1-\frac{\langle\partial\log\sigma_t^2(\boldsymbol{\theta}),\mathbf{u}\rangle}{2\sqrt{n}}\right)^k\right.\\ &\left.- \left(1-\frac{k\langle\partial\log\sigma_t^2(\boldsymbol{\theta}),\mathbf{u}\rangle}{2\sqrt{n}}\right)\right| = o_P(1).\end{aligned}$$

These are proven in Lemma 3.7. This completes the proof of Theorem 1.1. $\square$

As a consequence of Theorem 1.1, we can obtain the consistency result $\hat{\mu}_k = \hat{T}_n^{(k)}(1)/n \to \mu_k$ for $k\geq 2$, where $\hat{T}_n^{(k)}$ is given in (1.15). First, for any $1\leq i\leq k$, Theorem 1.1 implies

$$(3.8)\qquad \frac{1}{n}\sum_{t=1}^{n}\hat{\varepsilon}_t^i = \frac{1}{n}\sum_{t=1}^{n}\varepsilon_t^i - \frac{i\mu_i}{2\sqrt{n}}\langle\psi(\boldsymbol{\theta}),\sqrt{n}(\hat{\boldsymbol{\theta}}_n-\boldsymbol{\theta})\rangle + o_P\left(\frac{1}{\sqrt{n}}\right).$$

In particular, since $\mu_1 = 0$, we obtain

$$(3.9)\qquad |\bar{\hat{\varepsilon}} - \bar{\varepsilon}| = o_P\left(\frac{1}{\sqrt{n}}\right).$$

Since $\hat{\boldsymbol{\theta}}_n$ is $\sqrt{n}$ consistent, then for $k\geq 2$ and $E|\varepsilon_0|^k < \infty$, equation (3.8) implies that

$$(3.10)\qquad \left|\frac{\hat{T}_n^{(k)}(1)}{n} - \frac{T_n^{(k)}(1)}{n}\right| = O_P\left(\frac{1}{\sqrt{n}}\right).$$

Proof of Theorem 1.3. Equation (3.5) implies that, for any $b > 0$ and an integer $k \geq 1$,

$$(3.11) \quad \sup_{|\mathbf{u}| \leq b} \sum_{t=1}^{n} \left| \frac{\varepsilon_t^k(\mathbf{u}) - \varepsilon_t^k}{\sqrt{n}} + \frac{k\varepsilon_t^k}{2n} \langle \partial \log \sigma_t^2(\boldsymbol{\theta}), \mathbf{u} \rangle \right| = o_P(1).$$

By the inequality $||a| - |b|| \leq |a - b|$, (3.11) yields

$$\sup_{|\mathbf{u}| \leq b} \left| \frac{1}{\sqrt{n}} \sum_{t=1}^{n} |\varepsilon_t^k(\mathbf{u})| - \varepsilon_t^k| - \frac{k}{2n} \sum_{t=1}^{n} |\varepsilon_t^k \langle \partial \log \sigma_t^2(\boldsymbol{\theta}), \mathbf{u} \rangle| \right| = o_P(1).$$

Using ergodicity, we have, for each $\mathbf{u} \in \mathbb{R}^{p+q+1}$,

$$\psi_k^{(n)}(\mathbf{u}) = \frac{1}{n} \sum_{t=1}^{n} |\varepsilon_t^k \langle \partial \log \sigma_t^2(\boldsymbol{\theta}), \mathbf{u} \rangle| \to \psi_k(\mathbf{u}) \qquad \text{a.s.}$$

In Remark 1.3 it is argued that $\psi_k(\mathbf{u})$ is Lipschitz. With this method we obtain that $\psi_k^{(n)}(\mathbf{u})$ is also Lipschitz a.s. uniformly in $n$. Hence, one can obtain

$$\sup_{|\mathbf{u}| \leq b} |\psi_k^{(n)}(\mathbf{u}) - \psi_k(\mathbf{u})| = o_P(1).$$

Thus,

$$\sup_{|\mathbf{u}| \leq b} \left| \frac{1}{\sqrt{n}} \sum_{t=1}^{n} |\varepsilon_t^k(\mathbf{u})| - \varepsilon_t^k| - \frac{k}{2} \psi_k(\mathbf{u}) \right| = o_P(1).$$

Theorem 1.3 now follows standard arguments. This completes the proof of Theorem 1.3. $\square$

Proof of Theorem 1.4. When $k = 1$, we have

$$\frac{1}{\sqrt{n}} \left| \frac{\hat{T}_n^{(1)}(u)}{\hat{\sigma}_{(n)}} - \frac{T_n^{(1)}(u)}{\sigma_{(n)}} \right| \leq \frac{1}{\sqrt{n}} \frac{|\hat{T}_n^{(1)}(u) - T_n^{(1)}(u)|}{\hat{\sigma}_{(n)}} + \frac{|T_n^{(1)}(u)|}{\sqrt{n}} \left| \frac{1}{\hat{\sigma}_{(n)}} - \frac{1}{\sigma_{(n)}} \right|.$$

Since $E\varepsilon_0^2 < \infty$, the invariance principle for i.i.d. partial sums and (3.10) imply that

$$\sup_{0 \leq u \leq 1} \frac{|T_n^{(1)}(u)|}{\sqrt{n}} \left| \frac{1}{\hat{\sigma}_{(n)}} - \frac{1}{\sigma_{(n)}} \right| = o_P(1).$$

Hence, we can prove Theorem 1.4 for the case $k = 1$ if

$$\sup_{0 \leq u \leq 1} \frac{|\hat{T}_n^{(1)}(u) - T_n^{(1)}(u)|}{\sqrt{n}} = o_P(1),$$

which follows immediately by

$$\sup_{0 \le u \le 1} |\hat{T}_n^{(1)}(u) - T_n^{(1)}(u)|$$

$$\le \sup_{0 \le u \le 1} |(\hat{S}_n^{(1)}(u) - u\hat{S}_n^{(1)}(1)) - (S_n^{(1)}(u) - uS_n^{(1)}(1))| + |\hat{\bar{\varepsilon}} - \bar{\varepsilon}|$$

and by Theorem 1.2 and (3.9).

Next we consider the case $k \ge 2$. Let

$$\hat{L}_n(u) = \sqrt{n}\left(\frac{1}{n}\sum_{t=1}^{[nu]} \hat{T}_n^{(k)}(u) - u\lambda_k \hat{\sigma}_{(n)}^k\right)$$

and

$$L_n(u) = \sqrt{n}\left(\frac{1}{n}\sum_{t=1}^{[nu]} T_n^{(k)}(u) - u\lambda_k \sigma_{(n)}^k\right).$$

Then

$$\sup_{0 \le u \le 1} \frac{1}{\sqrt{n}}\left|\frac{\hat{T}_n^{(k)}(u)}{\hat{\sigma}_{(n)}^k} - \frac{T_n^{(k)}(u)}{\sigma_{(n)}^k}\right|$$

$$\le \frac{\sup_{0 \le u \le 1} |\hat{L}_n(u) - L_n(u)|}{\hat{\sigma}_{(n)}^k} + \sup_{0 \le u \le 1} |L_n(u)|\left|\frac{1}{\hat{\sigma}_{(n)}^k} - \frac{1}{\sigma_{(n)}^k}\right|.$$

Notice that (3.10) implies

$$\left|\frac{1}{\hat{\sigma}_{(n)}^k} - \frac{1}{\sigma_{(n)}^k}\right| = O_P\left(\frac{1}{\sqrt{n}}\right).$$

Therefore, we can prove Theorem 1.4 if

$$(3.12) \qquad \sup_{0 \le u \le 1} |\hat{L}_n(u) - L_n(u)| = o_P(1)$$

and

$$(3.13) \qquad \sup_{0 \le u \le 1} \frac{|L_n(u)|}{\sqrt{n}} = o_P(1).$$

By the facts that $\bar{\varepsilon} = O_P(1/\sqrt{n})$ and $\sigma_{(n)}^2 = \mu_2 + o_P(1)$, (3.13) holds if

$$\sup_{0 \le u \le 1}\left|\frac{1}{n}\sum_{t=1}^{[nu]} \varepsilon_t^k - u\mu_k\right| = o_P(1),$$

which is true by Lemma 3.6.

To prove (3.12), we need a finer representation of $\hat{\sigma}_{(n)}^k$ in terms of $\sigma_{(n)}^k$. By (3.8) and (3.9) we obtain

$$\hat{\sigma}_{(n)}^2 = \sigma_{(n)}^2 - \frac{\mu_2}{\sqrt{n}}\langle\psi(\boldsymbol{\theta}), \sqrt{n}(\hat{\boldsymbol{\theta}}_n - \boldsymbol{\theta})\rangle + o_P\left(\frac{1}{\sqrt{n}}\right).$$

By a first-order Taylor approximation with remainder we obtain

$$\hat{\sigma}_{(n)}^k = \left(\sigma_{(n)}^2 - \frac{\mu_2}{\sqrt{n}}\langle\psi(\boldsymbol{\theta}), \sqrt{n}(\hat{\boldsymbol{\theta}}_n - \boldsymbol{\theta})\rangle + o_P\left(\frac{1}{\sqrt{n}}\right)\right)^{k/2}$$

$$(3.14) \qquad = \left(\sigma_{(n)}^2 - \frac{\mu_2}{\sqrt{n}}\langle\psi(\boldsymbol{\theta}), \sqrt{n}(\hat{\boldsymbol{\theta}}_n - \boldsymbol{\theta})\rangle\right)^{k/2} + o_P\left(\frac{1}{\sqrt{n}}\right)$$

$$= \sigma_{(n)}^k - \frac{k\sigma_{(n)}^{k-2}}{2\sqrt{n}}\mu_2\langle\psi(\boldsymbol{\theta}), \sqrt{n}(\hat{\boldsymbol{\theta}}_n - \boldsymbol{\theta})\rangle + o_P\left(\frac{1}{\sqrt{n}}\right).$$

On the other hand, by (3.8) and (3.9), and the facts that $\bar{\hat{\varepsilon}} = O_P(1/\sqrt{n})$ and $\bar{\varepsilon} = O_P(1/\sqrt{n})$, we obtain

$$\frac{1}{\sqrt{n}}\hat{T}_n^{(k)}(u) = \frac{1}{\sqrt{n}}\sum_{t=1}^{[nu]}\hat{\varepsilon}_t^k - \frac{k\bar{\hat{\varepsilon}}}{\sqrt{n}}\sum_{t=1}^{[nu]}\hat{\varepsilon}_t^{k-1} + o_P(1)$$

$$= \frac{1}{\sqrt{n}}\sum_{t=1}^{[nu]}\varepsilon_t^k - \frac{ku\mu_k}{2}\langle\psi(\boldsymbol{\theta}), \sqrt{n}(\hat{\boldsymbol{\theta}}_n - \boldsymbol{\theta})\rangle - \frac{k\bar{\varepsilon}}{\sqrt{n}}\sum_{t=1}^{[nu]}\varepsilon_t^{k-1} + o_P(1)$$

$$= \frac{1}{\sqrt{n}}T_n^{(k)}(u) - \frac{ku\mu_k}{2}\langle\psi(\boldsymbol{\theta}), \sqrt{n}(\hat{\boldsymbol{\theta}}_n - \boldsymbol{\theta})\rangle + o_P(1)$$

uniformly in $0 \leq u \leq 1$. Substituting the above expression and (3.14) into (3.12) yields

$$\hat{L}_n(u) = \frac{1}{\sqrt{n}}T_n^{(k)}(u) - \frac{ku\mu_k}{2}\langle\psi(\boldsymbol{\theta}), \sqrt{n}(\hat{\boldsymbol{\theta}}_n - \boldsymbol{\theta})\rangle$$

$$- u\lambda_k\left(\sqrt{n}\sigma_{(n)}^k - \frac{k\mu_2^{k/2}}{2}\langle\psi(\boldsymbol{\theta}), \sqrt{n}(\hat{\boldsymbol{\theta}}_n - \boldsymbol{\theta})\rangle\right) + o_P(1)$$

$$= L(u) + o_P(1)$$

uniformly in $0 \leq u \leq 1$. This concludes the proof of (3.12) and, hence, Theorem 1.4. $\square$

The remainder of this section gives the various lemmas needed in the proofs above, plus the proof of Theorem 1.5.

By the mean value theorem for the multivariate function $\sigma_t^2(\mathbf{u})$ there is a $\boldsymbol{\zeta}$ satisfying $|\boldsymbol{\zeta} - \boldsymbol{\theta}| \leq |\mathbf{u}|/\sqrt{n}$ so that (3.3) yields

$$|g_t(\mathbf{u})| = |\mathbf{u}|\frac{|\partial\sigma_t^2(\boldsymbol{\zeta})|}{\sigma_t^2(\boldsymbol{\theta})} = |\mathbf{u}|\ \frac{|\partial\sigma_t^2(\boldsymbol{\zeta})|}{\sigma_t^2(\boldsymbol{\zeta})}\frac{\sigma_t^2(\boldsymbol{\xi})}{\sigma_t^2(\boldsymbol{\theta})}.$$

Also, by a second-order term Taylor expansion there exists $\boldsymbol{\xi}$ satisfying $|\boldsymbol{\xi} - \boldsymbol{\theta}| \leq |\mathbf{u}|/\sqrt{n}$ so that

$$g_t(\mathbf{u}) = \langle \partial \log \sigma_t^2(\boldsymbol{\theta}), \mathbf{u} \rangle + \frac{1}{2\sqrt{n}} \mathbf{u} \frac{\partial^2 \sigma_t^2(\boldsymbol{\xi})}{\sigma_t^2(\boldsymbol{\theta})} \mathbf{u}^\tau,$$

where $\mathbf{u}^\tau$ is the transpose of the vector $\mathbf{u}$, and $\partial^2 \sigma_t^2(\mathbf{u})$ is the matrix of the second-order partial derivatives of $\sigma_t^2(\mathbf{u})$ (the Hessian matrix). Therefore,

$$\sup_{|\mathbf{u}| \leq b} |g_t(\mathbf{u}) - \langle \partial \log \sigma_t^2(\boldsymbol{\theta}), \mathbf{u} \rangle| \leq \frac{b^2}{2\sqrt{n}} \sup_{|\boldsymbol{\xi} - \boldsymbol{\theta}| \leq b/\sqrt{n}} \frac{|\partial^2 \sigma_t^2(\boldsymbol{\xi})|}{\sigma_t^2(\boldsymbol{\xi})} \frac{\sigma_t^2(\boldsymbol{\xi})}{\sigma_t^2(\boldsymbol{\theta})}.$$

Thus, to show that $g_t(\mathbf{u})$ has a finite moment and can be approximated by $\langle \partial \log \sigma_t^2(\boldsymbol{\theta}), \mathbf{u} \rangle$ in the neighborhood of $|\mathbf{u}| \leq b$ for some $b > 0$, we need the following two lemmas.

LEMMA 3.1. *If* (1.6) *and* (1.4) *hold, then* $E|\varepsilon_0|^\delta < \infty$ *for some* $\delta > 0$ *implies that*

$$E\left( \sup_{\mathbf{u} \in \boldsymbol{\Theta}} |\partial \log \sigma_t^2(\mathbf{u})| \right)^{\kappa^*} = E\left( \sup_{\mathbf{u} \in \boldsymbol{\Theta}} \frac{|\partial \sigma_t^2(\mathbf{u})|}{\sigma_t^2(\mathbf{u})} \right)^{\kappa^*} < \infty$$

*and*

$$E\left( \sup_{\mathbf{u} \in \boldsymbol{\Theta}} \frac{|\partial^2 \sigma_t^2(\mathbf{u})|}{\sigma_t^2(\mathbf{u})} \right)^{\kappa^*} < \infty$$

*for any* $\kappa^* > 0$.

PROOF. See the proof of Lemma 5.6 of [3]. $\square$

LEMMA 3.2. *If* (1.6) *and* (1.14) *hold, then* $E|\varepsilon_0|^\delta < \infty$ *for some* $\delta > 0$ *implies that, for any* $b > 0$ *and* $\kappa^* > 0$, *there exists an integer* $N$ *such that*

$$\sup_{n \geq N} E\left( \sup_{|\mathbf{u}| \leq b} \frac{\sigma_t^2(\boldsymbol{\theta} + n^{-1/2}\mathbf{u})}{\sigma_t^2(\boldsymbol{\theta})} \right)^{\kappa^*} < \infty.$$

PROOF. By Lemma 3.1 of [3], when $n$ is large enough (so that $\boldsymbol{\theta} + n^{-1/2}\mathbf{u} \in \boldsymbol{\Theta}$), we have

$$0 \leq c_i(\boldsymbol{\theta} + n^{-1/2}\mathbf{u}) \leq C_1 \left( \max_{1 \leq j \leq q} \frac{\beta_j + n^{-1/2}|t_j|}{\beta_j} \right)^i c_i(\boldsymbol{\theta}) \leq C_1 \rho_n^i c_i(\boldsymbol{\theta}),$$
$$0 \leq i < \infty,$$

where $C_1$ is a constant and $1 < \rho_n = 1 + n^{-1/2}b/\underline{u}$, where $\underline{u}$ is defined above (1.6). Thus, Lemma 3.2 will be proven if we can show that

$$E\left( \frac{\sum_{i=1}^{\infty} \rho_N^i c_i(\boldsymbol{\theta}) X_{t-i}^2}{1 + \sum_{i=1}^{\infty} c_i(\boldsymbol{\theta}) X_{t-i}^2} \right)^{\kappa^*} < \infty.$$

Since $\rho_N$ can be close enough to 1 if $N$ is large enough, the above inequality follows from the same proof of Lemma 3.7 of [2]. For completeness, we give a detailed proof here.

By Lemma 3.1 of [3], there are constants $C_2$ and $0 < \rho < 1$ such that

$$|c_i(\mathbf{u})| \leq C_2 \rho^i \qquad \text{for all } \mathbf{u} \in \boldsymbol{\Theta} \text{ and all } i.$$

Then for any $M \geq 1$, we have

$$\frac{\sum_{i=1}^{\infty} \rho_N^i c_i(\boldsymbol{\theta}) X_{t-i}^2}{1 + \sum_{i=1}^{\infty} c_i(\boldsymbol{\theta}) X_{t-i}^2} \leq \rho_N^M + \sum_{i=M+1}^{\infty} \rho_N^i c_i(\boldsymbol{\theta}) X_{t-i}^2$$

$$\leq \rho_N^M + C_2 \sum_{i=M+1}^{\infty} (\rho_N \rho)^i X_{t-i}^2.$$

By Lemma 2.3 of [3],

(3.15)        there exists $\delta^* > 0$ such that $E|X_0|^{\delta^*} < \infty$.

Notice that, for $N$ sufficiently large, $\rho_N \rho < 1$. By Markov's inequality, we have, for $x > \rho_N^2$,

$$P\left\{ \sum_{i=M+1}^{\infty} (\rho_N \rho)^i X_{t-i}^2 > x/2 \right\}$$

$$\leq \sum_{i=M+1}^{\infty} P\{ X_{t-i}^2 > (x/2)(\rho_N \rho)^{-i}(1 - (\rho_N \rho)^{1/2})(\rho_N \rho)^{i/2} \}$$

$$= \sum_{i=M+1}^{\infty} P\{ |X_{t-i}|^{\delta^*} > (x/2)^{\delta^*/2}(1 - (\rho_N \rho)^{1/2})^{\delta^*/2}(\rho_N \rho)^{-i\delta^*/4} \}$$

$$\leq E|X_0|^{\delta^*}(x/2)^{-\delta^*/2}(1 - (\rho_N \rho)^{1/2})^{-\delta^*/2}(1 - (\rho_N \rho)^{\delta^*/4})^{-1}(\rho_N \rho)^{M\delta^*/4}.$$

Choosing $M = \log(C_2 x/2)/\log \rho_N$, we have, for any $\kappa^* > 0$,

$$P\left\{ \frac{\sum_{i=1}^{\infty} \rho_N^i c_i(\boldsymbol{\theta}) X_{t-i}^2}{1 + \sum_{i=1}^{\infty} c_i(\boldsymbol{\theta}) X_{t-i}^2} > C_2 x \right\}$$

$$\leq P\left\{ \sum_{i=M+1}^{\infty} (\rho_N \rho)^i X_{t-i}^2 > x/2 \right\}$$

$$\leq C_3 \exp(-(\delta^*/4)(1 + \log \rho^{-1}/\log \rho_N) \log(x/2))$$

$$\leq C_4 x^{-2\kappa^*}$$

if $\rho_N > 1$ is close enough to 1 (when $N$ is sufficiently large), and where $C_3$, $C_4$ are constants that may depend on $N$.   $\square$

Using Lemmas 3.1 and 3.2 and Hölder's inequality, we immediate arrive at the following result.

LEMMA 3.3. *If* (1.6) *and* (1.14) *hold, then* $E|\varepsilon_0|^\delta < \infty$ *for some* $\delta > 0$ *implies that, for any* $b > 0$ *and* $\kappa^* > 0$, *there exists an integer* $N$ *such that*

$$\sup_{n \geq N} E\Big(\sup_{|\mathbf{u}| \leq b} |g_t(\mathbf{u})|\Big)^{\kappa^*} < \infty$$

*and*

$$\sup_{n \geq N} E\Big(\sup_{|\mathbf{u}| \leq b} \sqrt{n}|g_t(\mathbf{u}) - \langle \partial \log \sigma_t^2(\boldsymbol{\theta}), \mathbf{u}\rangle|\Big)^{\kappa^*} < \infty.$$

LEMMA 3.4. *Suppose that* (1.6) *and* (1.14) *hold and* $E|\varepsilon_0|^\delta < \infty$ *for some* $\delta > 0$. *Then for any* $b > 0$,

$$\max_{1 \leq t \leq n} \sup_{|\mathbf{u}| \leq b} \Big| \frac{1}{\sqrt{1 + n^{-1/2}g_t(\mathbf{u})}} - \Big(1 - \frac{\langle \partial \log \sigma_t^2(\boldsymbol{\theta}), \mathbf{u}\rangle}{2\sqrt{n}}\Big)\Big| = o_P\Big(\frac{1}{\sqrt{n}}\Big).$$

PROOF. First we state a well-known result that if $\{Y_n, n \geq 0\}$ is a sequence of identically distributed r.v.'s with $E|Y_0|^{\kappa^*} < \infty$ for some $\kappa^* > 0$, then

(3.16) $$\max_{1 \leq t \leq n} |Y_t| = o_P(n^{1/\kappa^*}).$$

Thus, by Lemma 3.3 for $\kappa^* > 2$, and noting that by construction $g_t(\mathbf{u})$ has the same marginal distribution for each $t$, we have

$$n^{-1/2} \max_{1 \leq t \leq n} \sup_{|\mathbf{u}| \leq b} |g_t(\mathbf{u})| = o_P(n^{1/\kappa^* - 1/2}).$$

This, together with the inequality that $|1/\sqrt{1 + x} - 1 + x/2| \leq 3x^2$ for $|x| \leq 1/2$, implies that

$$\max_{1 \leq t \leq n} \sup_{|\mathbf{u}| \leq b} \Big| \frac{1}{\sqrt{1 + n^{-1/2}g_t(\mathbf{u})}} - \Big(1 - \frac{g_t(\mathbf{u})}{2\sqrt{n}}\Big)\Big| = (o_P(n^{1/\kappa^* - 1/2}))^2 = o_P(n^{2/\kappa^* - 1}).$$

By Lemma 3.3 and (3.16),

$$\max_{1 \leq t \leq n} \sup_{|\mathbf{u}| \leq b} \sqrt{n}|g_t(\mathbf{u}) - \langle \partial \log \sigma_t^2(\boldsymbol{\theta}), \mathbf{u}\rangle| = o_P(n^{1/\kappa^*}).$$

Choosing $\kappa^* > 4$ completes the proof of Lemma 3.4. □

LEMMA 3.5. *Suppose that* (1.6) *and* (1.14) *hold and that* $E(|\varepsilon_0|^\delta) < \infty$ *for some* $\delta > 0$. *Then for any integers* $k, \ell \geq 1$,

$$\sum_{t=1}^n |\tilde{\varepsilon}_t|^k |\tilde{\sigma}_t - \hat{\sigma}_t|^\ell = O_P(1).$$

PROOF. By (1.6) and (1.5), for any small $\eta > 0$ there exist a set with probability $> 1 - \eta$ and a constant $C_5$ such that

$$\tilde{\sigma}_t \geq C_5 \quad \text{and} \quad \hat{\sigma}_t \geq C_5 \qquad \text{for all } t \geq 1.$$

Thus, by (1.7), (1.8) and (1.9), there is a constant $C_6$ such that

$$|\tilde{\varepsilon}_t|^k |\tilde{\sigma}_t - \hat{\sigma}_t|^\ell \leq (2C_5)^{-\ell} |\varepsilon_t|^k \sup_{u \in \Theta} \left( \frac{\sigma_t(\theta)}{\sigma_t(u)} \right)^k \left( \sup_{u \in \Theta} \sum_{i=t+1}^\infty c_i(u) X_{t-i}^2 \right)^\ell$$

$$\leq (2C_5)^{-\ell} |\varepsilon_t|^k \sup_{u \in \Theta} \left( \frac{\sigma_t^2(\theta)}{\sigma_t^2(u)} \right)^{k/2} \left( C_6 \sum_{i=t+1}^\infty \rho^i X_{t-i}^2 \right)^\ell.$$

By Lemma 5.1 of [3], taking $0 < \nu = \delta/4 < \delta/2$ in their lemma,

$$E \left( \sup_{u \in \Theta} \frac{\sigma_t^2(\theta)}{\sigma_t^2(u)} \right)^{\delta/4} < \infty.$$

By (3.15), and taking $\delta^*/2 \leq 1$,

$$E \left[ \left( \sum_{i=t+1}^\infty \rho^i X_{t-i}^2 \right)^{\delta^*/2} \right] \leq \sum_{i=t+1}^\infty \rho^{i\delta^*/2} E |X_{t-i}|^{\delta^*}$$

$$\leq E|X_0|^{\delta^*} \frac{\rho^{t\delta^*/2}}{1 - \rho^{\delta^*/2}}.$$

Hence, by Hölder's inequality

$$E \left( \sum_{t=1}^n |\tilde{\varepsilon}_t|^k |\tilde{\sigma}_t - \hat{\sigma}_t|^\ell \right)^{\delta^{**}} \leq E \left( \sum_{t=1}^\infty |\tilde{\varepsilon}_t|^k |\tilde{\sigma}_t - \hat{\sigma}_t|^\ell \right)^{\delta^{**}} < \infty$$

for sufficiently small $\delta^{**} > 0$. Thus, Lemma 3.5 is now proven. $\quad\square$

LEMMA 3.6. *Let* $Y_t = h(\varepsilon_t, \varepsilon_{t-1}, \dots)$ *and suppose that* $E|Y_0| < \infty$. *Then*

$$\sup_{0 \leq u \leq 1} \left| \frac{1}{n} \sum_{t=1}^{[nu]} Y_t - u E Y_0 \right| = o_P(1).$$

PROOF. For any $0 < \xi < 1$ and for large $n$,

$$\sup_{0 \le u \le 1} \left| \frac{1}{n} \sum_{t=1}^{[nu]} Y_t - uEY_0 \right| \le \frac{1}{n} \sum_{t=1}^{[n\xi]} |Y_t| + \xi E|Y_0|$$

$$+ \sup_{\xi \le u \le 1} \left| \frac{1}{n} \sum_{t=1}^{[nu]} Y_t - uEY_0 \right| \le \left| \frac{1}{n} \sum_{t=1}^{[n\xi]} (|Y_t| - E|Y_0|) \right|$$

$$+ 3\xi E|Y_0| + \sup_{j \ge [n\xi]} \left| \frac{1}{j} \sum_{t=1}^{j} (Y_t - EY_0) \right|.$$

Since $Y_t = h(\varepsilon_t, \varepsilon_{t-1}, \ldots)$, by Theorem 3.5.8 of [22] $\{Y_t\}$ is stationary and ergodic. Thus, as $n \to \infty$,

$$\frac{1}{[n\xi]} \sum_{t=1}^{[n\xi]} (|Y_t| - E|Y_0|) = o_P(1)$$

and

$$\sup_{j \ge [n\xi]} \left| \frac{1}{j} \sum_{t=1}^{j} (Y_t - EY_0) \right| = o_P(1). \qquad \square$$

LEMMA 3.7. *Suppose that* (1.6) *and* (1.14) *hold. Then, for any $b > 0$ and an integer $k \ge 1$, $E|\varepsilon_0|^k < \infty$ implies:*

(i)
$$\sup_{|\mathbf{u}| \le b} \sum_{t=1}^{n} |\varepsilon_t|^k \left| \left( 1 - \frac{\langle \partial \log \sigma_t^2(\boldsymbol{\theta}), \mathbf{u} \rangle}{2\sqrt{n}} \right)^k - \left( 1 - \frac{k \langle \partial \log \sigma_t^2(\boldsymbol{\theta}), \mathbf{u} \rangle}{2\sqrt{n}} \right) \right| = O_P(1)$$

*and*

(ii)
$$\sup_{|\mathbf{u}| \le b} \frac{1}{\sqrt{n}} \sum_{t=1}^{n} |\varepsilon_t|^k \left| \left( \frac{1}{\sqrt{1 + n^{-1/2} g_t(\mathbf{u})}} \right)^k - \left( 1 - \frac{\langle \partial \log \sigma_t^2(\boldsymbol{\theta}), \mathbf{u} \rangle}{2\sqrt{n}} \right)^k \right| = o_P(1).$$

PROOF. Part (i), case $k = 1$ is trivial. Part (ii), case $k = 1$ follows directly from Lemma 3.4.

It is easy to see from Lemma 3.1 that

$$E \sup_{|\mathbf{u}| \le b} |\langle \partial \log \sigma_t^2(\boldsymbol{\theta}), \mathbf{u} \rangle|^k < \infty.$$

Now consider $k \ge 2$. Thus, using the fact that $\varepsilon_t$ and $\partial \log \sigma_t^2(\boldsymbol{\theta})$ are independent, we have, for $2 \le i \le k$,

$$E \left( \sup_{|\mathbf{u}| \le b} \frac{1}{n} \sum_{t=1}^{n} |\varepsilon_t|^k |\langle \partial \log \sigma_t^2(\boldsymbol{\theta}), \mathbf{u} \rangle|^i \right)$$

$$\leq \frac{1}{n} \sum_{t=1}^{n} E|\varepsilon_t|^k E \sup_{|\mathbf{u}| \leq b} |\langle \partial \log \sigma_t^2(\boldsymbol{\theta}), \mathbf{u} \rangle|^i$$

$$= O(1).$$

This, together with the binomial formula, implies Lemma 3.7(i).

Now consider part (ii) and $k \geq 2$. By Lemma 3.4 we have

$$\frac{1}{\sqrt{1 + n^{-1/2} g_t(\mathbf{u})}} = 1 - \frac{\langle \partial \log \sigma_t^2(\boldsymbol{\theta}), \mathbf{u} \rangle}{2\sqrt{n}} + o_P\left(\frac{1}{\sqrt{n}}\right)$$

uniformly in $|\mathbf{u}| \leq b$ and $1 \leq t \leq n$. Thus, using the inequality

$$|(x + \Delta)^k - x^k| \leq k 2^{k-1} |\Delta| (|x|^{k-1} + |\Delta|^{k-1}),$$

we have

$$\left| \left( \frac{1}{\sqrt{1 + n^{-1/2} g_t(\mathbf{u})}} \right)^k - \left( 1 - \frac{\langle \partial \log \sigma_t^2(\boldsymbol{\theta}), \mathbf{u} \rangle}{2\sqrt{n}} \right)^k \right|$$

$$\leq o_P\left(\frac{1}{\sqrt{n}}\right) \left( \left| 1 - \frac{\langle \partial \log \sigma_t^2(\boldsymbol{\theta}), \mathbf{u} \rangle}{2\sqrt{n}} \right|^{k-1} + o_P\left(\frac{1}{\sqrt{n}}\right) \right)$$

uniformly in $|\mathbf{u}| \leq b$ and $1 \leq t \leq n$. Finally, using the fact that $\varepsilon_t$ and $\partial \log \sigma_t^2(\boldsymbol{\theta})$ are independent and $E \sup_{|\mathbf{u}| \leq b} |\langle \partial \log \sigma_t^2(\boldsymbol{\theta}), \mathbf{u} \rangle|^k < \infty$, the proof of Lemma 3.7(ii) is immediate.  □

LEMMA 3.8.  *If $E\varepsilon_0^{2k} < \infty$ for an integer $k \geq 2$, then*

$$\left\{ \frac{1}{\sqrt{n}} \left( \frac{T_n^{(k)}(u)}{\sigma_{(n)}^k} - n u \lambda_k \right), 0 \leq u \leq 1 \right\}$$

*converges weakly to the Gaussian process $\{B^{(k)}(u), 0 \leq u \leq 1\}$ with covariance defined by* (1.16).

PROOF.  Using the standard GARCH scaling assumption $\mu_2 = 1$, we have $\lambda_k = \mu_k$. Similar to the proof of Lemma 3.7, by Lemma 3.6 and $\bar{\varepsilon} = O_P(1/\sqrt{n})$ we have

$$\frac{1}{n} \sum_{t=1}^{[nu]} (\varepsilon_t - \bar{\varepsilon})^k = \frac{1}{n} \sum_{t=1}^{[nu]} \varepsilon_t^k - \frac{k}{n} \sum_{t=1}^{[nu]} \varepsilon_t^{k-1} \bar{\varepsilon} + o_P\left(\frac{1}{\sqrt{n}}\right)$$

$$= \frac{1}{n} \sum_{t=1}^{[nu]} (\varepsilon_t^k - \mu_k) - u k \mu_{k-1} \bar{\varepsilon} + u \mu_k + o_P\left(\frac{1}{\sqrt{n}}\right)$$

uniformly in $0 \le u \le 1$. On the other hand, by $\bar{\varepsilon}^2 = O_P(1/n)$ we have

$$
\left( \frac{1}{n} \sum_{t=1}^{n} (\varepsilon_t - \bar{\varepsilon})^2 \right)^{k/2} = \left( 1 + \frac{1}{n} \sum_{t=1}^{n} (\varepsilon_t^2 - 1) \right)^{k/2} + o_P\left( \frac{1}{\sqrt{n}} \right)
$$

$$
= 1 + \frac{k}{2n} \sum_{t=1}^{n} (\varepsilon_t^2 - 1) + o_P\left( \frac{1}{\sqrt{n}} \right).
$$

Therefore,

$$
\frac{1}{\sqrt{n}} \left( \frac{T_n^{(k)}(u)}{\sigma_{(n)}^k} - nu\lambda_k \right)
$$

$$
= \frac{1}{\sigma_{(n)}^k \sqrt{n}} \left( \sum_{t=1}^{[nu]} (\varepsilon_t^k - \mu_k) - \frac{uk}{2} \sum_{t=1}^{n} (\mu_k(\varepsilon_t^2 - 1) + 2\mu_{k-1}\varepsilon_t) \right) + o_P(1)
$$

$$
= \frac{1}{\sigma_{(n)}^k} M_n^{(k)}(u) + o_P(1)
$$

uniformly in $0 \le u \le 1$. Since $\hat{\sigma}_{(n)}^k \to 1$ in probability, we can prove Lemma 3.8 if $\{M_n^{(k)}(u), 0 \le u \le 1\}$ converges weakly to the Gaussian process $\{B^{(k)}(u), 0 \le u \le 1\}$. According to the invariance principle for i.i.d. partial sums, $\{n^{-1/2} \sum_{t=1}^{[nu]} (\varepsilon_t^k - \mu_k), 0 \le u \le 1\}$ and $\{n^{-1/2}u \sum_{t=1}^{n} (\mu_k(\varepsilon_t^2 - 1) + 2\mu_{k-1}\varepsilon_t), 0 \le u \le 1\}$ are tight (in fact, they converge weakly to Gaussian processes) and, hence, so is the process $\{M_n^{(k)}(u), 0 \le u \le 1\}$. Using the form of $M_n^{(k)}(u)$ and relatively easy but lengthy covariance computations, we obtain (1.16). These straightforward details are omitted. This completes the proof of Lemma 3.8. □

PROOF OF THEOREM 1.5. Using the form of $M_n^{(k)}(u)$ from the proof of Lemma 3.8, we just need to show that, for any $0 \le u, v \le 1$, as $n \to \infty$,

$$
EM_n^{(k)}(u) M_n^{(k+1)}(v) \to 0.
$$

This computation is straightforward but lengthy. The details are omitted. □

**Acknowledgments.** The authors thank three anonymous referees for their detailed comments and suggestions. These have led to improved and streamlined proofs and some remarks in this paper.

Department of Statistical
and Actuarial Sciences
University of Western Ontario
London, Ontario
Canada N6A 5B7
e-mail: kulperger@uwo.ca